\documentclass[11pt,leqno]{article}
\usepackage[centertags]{amsmath}
\usepackage{amsfonts}
\usepackage{amsmath}
\usepackage{amssymb}
\usepackage{latexsym}
\usepackage{amsthm}
\usepackage{mathrsfs}
\usepackage{newlfont}
\usepackage[latin1]{inputenc}

\usepackage{epsfig}
\usepackage{pst-all}
\usepackage{graphicx}
\usepackage{pstricks}

\newtheorem{teor}{Theorem}[section]

\newtheorem{lema}{Lemma}[section]

\newtheorem{coro}{Corollary}[section]

\newtheorem{nota}{Remark}[section]

\newcommand\RR{{\mathbb R}}

\newcommand{\C}{\mathbb{C}^{N\times N}}

\newcommand{\Ld}{L^2([a,b],\mathbb{C}^{N\times N})}
\newcommand{\LdR}{L^2(\mathbb{R},\mathbb{C}^{N\times N})}
\newcommand{\LWd}{L_W^2([a,b],\mathbb{C}^{N\times N})}

\newcommand{\pro}{\par{\scshape Proof:} }

\numberwithin{equation}{section} 

\title{Some examples of matrix-valued orthogonal functions having a differential and an integral operator as eigenfunctions\thanks{The work of the author is partially supported by D.G.E.S, ref. MTM2009-12740-C03-02, Junta de Andaluc\'{i}a, grants FQM-229, FQM-481, P06-FQM-01738 and Subprograma de estancias de movilidad posdoctoral en el extranjero, MICINN, ref. -2008-0207.\newline \textbf{Mathematics Subject
Classification (2010)} 42C05, 34L40. \newline \textbf{Keywords} Matrix-valued Schr\"odinger operators, matrix-valued orthogonal polynomials, Fourier analysis}}

\author{Manuel D. de la Iglesia \\ \footnotesize Courant Institute of Mathematical Sciences. New York University \\ \footnotesize\ 251 Mercer Street, New York, NY 10012, U.S.A. mdi29@cims.nyu.edu}

\date{}

\begin{document}

\maketitle

\begin{abstract}
The aim of this paper is to show some examples of matrix-valued orthogonal functions on the real line which are simultaneously eigenfunctions of a second-order differential operator of Schr\"{o}dinger type and an integral operator of Fourier type. As a consequence we derive integral equations of these functions as well as other useful structural formulas. Some of these functions are plotted to show the relationship with the Hermite or wave functions.
\end{abstract}

\section{Introduction}\label{INTRO}

In this paper we will show examples of $N\times N$ matrix-valued orthogonal functions $(\Phi_n)_n$ which are simultaneously eigenfunctions of a second-order differential operator of Schr\"{o}dinger type, i.e.
\begin{equation}\label{Schrodgen}
    (\Phi_n\mathcal{D})(x)\doteq\Phi_n''(x)-\Phi_n(x)V(x)=\Gamma_n\Phi_n(x),\quad x\in\mathbb{R},
\end{equation}
where the matrix-valued potential $V(x)$ (independent of $n$) is a diagonal quadratic matrix polynomial (but not a scalar multiple of the identity), and an integral operator of Fourier type, i.e.
\begin{equation}\label{Inteqgen}
    (\Phi_n\mathcal{I})(x)\doteq\frac{1}{\sqrt{2\pi}}\int_{-\infty}^{\infty}\Phi_n(t)K(x,t)dt=\Lambda_n\Phi_n(x),\quad x\in\mathbb{R},
\end{equation}
where the matrix-valued kernel $K(x,t)$ is also diagonal (but not a scalar multiple of the identity) of the form $K(x,t)=e^{ixt}\widetilde{K}$ for some diagonal matrix $\widetilde{K}$. We will also show that it is possible to construct suitable families $(\Phi_n)_n$ such that the corresponding eigenvalues $\Gamma_n$ and $\Lambda_n$ are diagonal. Hence both differential and integral operators commute in the space of all matrix-valued functions spanned by $(\Phi_n)_n$.

Observe that the operators $\mathcal{D}$ and $\mathcal{I}$ appear on the \emph{right}. This denotes that the potential coefficient and the kernel are multiplied on the \emph{right}, respectively, while the eigenvalues appear on the \emph{left}. With this configuration it is straightforward to derive that both differential and integral operators commute in the space of all matrix-valued functions spanned by $(\Phi_n)_n$, i.e. for any $F=\sum_{n=0}^mC_n\Phi_n, m\geq0, C_n\in\C$,
$$
 F\mathcal{I}\mathcal{D}=\Lambda_nF\mathcal{D}=\Lambda_n\Gamma_nF=\Gamma_n\Lambda_nF=\Gamma_nF\mathcal{I}=F\mathcal{D}\mathcal{I},
$$
since the eigenvalues are diagonal matrices. In particular we will see in Section \ref{CON} that this commutativity property holds in $\LdR$ (see Section \ref{PRE} for definitions) for the examples we consider in this paper.

In the scalar situation the examples of simultaneously eigenfunctions of an integral and a differential operator are very special and usually are related to many areas of engineering, physics, chemistry and mathematics. Prominent examples are the Hermite or wave functions which satisfy the Schr\"{o}dinger equation \eqref{Schrodgen} with $V(x)=x^2$ and also are eigenfunctions of the Fourier transform (see \cite{DMc, T, Wi}). The Hermite functions play an important role in quantum mechanics to model the one dimensional quantum harmonic oscillator (see \cite{PW}), or in chemistry in vibration of molecules (see \cite{PW, WDC}), among other applications.

Also the prolate spheroidal wave functions in signal processing satisfy certain second-order differential equation and are eigenfunctions of the so-called \emph{sinc kernel}. This remarkable fact was discovered in a series of papers by the Bell Labs group in the early 1960s, see \cite{SP, LP1, LP2, S1, S2}. See also \cite{Dau}.

The commuting property and other similar properties are very closely related to the classical orthogonal polynomials of Hermite, Laguerre, Jacobi and Bessel. For instance, in \cite{G3}, the author produces naturally appearing global operators that happen to commute with properly chosen local operators.

\bigskip

In the last few years many examples of matrix-valued orthogonal polynomials  $(P_n)_n$ have appeared satisfying second-order differential equations of Sturm-Liouville type, i.e.
\begin{equation*}\label{sodeintro}
    P_n''(x)F_2(x)+P_n'(x)F_1(x)+P_n(x)F_0(x)=\Gamma_nP_n(x),
\end{equation*}
where $F_2, F_1$ and $F_0$ are matrix polynomials (which do not depend on $n$) of degrees less than or equal to 2, 1 and 0, respectively, and $\Gamma_n$ are Hermitian matrices (if the family $(P_n)_n$ is orthonormal). These examples are the matrix analogue of the classical families of Hermite, Laguerre and Jacobi polynomials. Since the classical families are intimately related with integral equations, it is natural to expect that these new families of matrix-valued orthogonal polynomials are related to matrix-valued integral equations. The main goal of this paper is to show some examples of simultaneously eigenfunctions of an integral and a differential operator in the matrix case.

For that purpose we will consider families of matrix-valued orthogonal functions (constructed from the matrix-valued orthogonal polynomials and a weight matrix) such that they are orthogonal in the corresponding function spaces (more on this in Section \ref{PRE}). This treatment will allow us to obtain, for the examples we study in this paper, matrix-valued Schr\"odinger operators of the form \eqref{Schrodgen} and integral equations of the form \eqref{Inteqgen}.

Schr\"odinger operators with Hermitian matrix-valued potentials (usually with some restrictions on the potential) are not new in the literature. They have been considered in \cite{CGe, CGHL, GKM, OMG} to study many analytical properties of these operators. Also matrix-valued functions are considered for some problems in wavelet theory to model matrix-valued signals, e.g. video images (see \cite{WS, XS}).


\bigskip

To give an idea of the results contained in this paper let us display here one of the examples we study in Section \ref{FIRST}. Consider a family of matrix-valued functions $(\Phi_n)_n$ defined in \eqref{Phi1}. Then they satisfy the following matrix-valued Schr\"odinger equation
\begin{equation*}\label{Schrodinger1intro}
\Phi_n''(x)-\Phi_n(x)(x^2I+2J)+((2n+1)I+2J)\Phi_n(x)=0,\quad x\in\mathbb{R},
\end{equation*}
where $J$ is the diagonal matrix defined in \eqref{JJ}. This differential equation (for the monic family) was already derived in Section 6.2 of \cite{DG3} (with non diagonal eigenvalue). The family we consider in \eqref{Phi1} will transform the eigenvalue into another convenient diagonal one. Then the family of matrix-valued orthogonal functions $(\Phi_n)_n$ satisfies the following integral equation
\begin{equation*}
\frac{1}{\sqrt{2\pi}}\int_{-\infty}^{\infty}\Phi_n(t)e^{ixt}e^{i\frac{\pi}{2}J}dt=(i)^ne^{i\frac{\pi}{2}J}\Phi_n(x),\quad x\in\mathbb{R},
\end{equation*}
where $e^{i\frac{\pi}{2}J}$ is the diagonal matrix \eqref{iJ}.

Therefore the kernel and the eigenvalue in \eqref{Inteqgen} are $K(x,t)=e^{ixt}e^{i\frac{\pi}{2}J}$ and $\Lambda_n=(i)^ne^{i\frac{\pi}{2}J}$, respectively. Observe that they are diagonal, but not necessarily scalars multiple of the identity. This result will have important consequences. In particular we will derive some structural formulas and integral equations of the corresponding matrix-valued orthogonal polynomials.

\bigskip

The paper is organized as follows: in Section \ref{PRE} we give some preliminaries. In Section \ref{MVOF} we derive the matrix-valued differential equation satisfied by the matrix-valued functions constructed from the matrix-valued orthogonal polynomials and a weight matrix. We will focus on the case when the leading coefficient $F_2$ is of the form $F_2(x)=f_2(x)I$, where $f_2(x)$ is a scalar polynomial of degree less than or equal to 2 (the case considered in \cite{DG1}). This result holds for any family supported in any interval $[a,b]$, $-\infty\leq a<b\leq\infty$.

In Sections \ref{FIRST} and \ref{SECOND} we study in full detail two examples supported in $\mathbb{R}$, which previously appeared in \cite{DG1} and derive differential and integral equations for the corresponding families of matrix-valued orthogonal functions. As a consequence we will derive \emph{real} (and complex) integral equations for the families of matrix-valued orthogonal polynomials, among other structural formulas of general size $N\times N$. We also study, for both examples, the special case of $N=2$. In this case we can write our (normalized) matrix-valued orthogonal functions $\Phi_n$ in terms of classical Hermite or wave functions. We will plot some of the entries of $\Phi_n\Phi_n^*$ for the first values of $n$, as well as derive some integrals involving these functions. Finally, Section \ref{CON} addresses the issue of the completeness of these functions, gives a summary of the results of the paper and the challenges that lie ahead.

\section{Preliminaries}\label{PRE}

An $N\times N$ matrix-valued function $F$ is a matrix with real-valued functions as entries of the form
$$
F(x)=\begin{pmatrix}
                f_{11}(x) & f_{12}(x) & \cdots & f_{1N}(x) \\
                f_{21}(x) & f_{22}(x) & \cdots & f_{2N}(x) \\
                \vdots & \vdots & \ddots & \vdots \\
                f_{N1}(x) & f_{N2}(x) & \cdots & f_{NN}(x) \\
              \end{pmatrix},\quad x\in[a,b].
$$
We will say that $F\in\Ld$ if
\begin{equation*}\label{Ld}
    \int_{a}^{b}F(x)F^*(x)dx<\infty,
\end{equation*}
where $F^*$ denotes the Hermitian conjugate of the matrix $F$ and $-\infty\leq a,b\leq\infty$. In the above definition we mean that the integral is finite entry by entry. It is easy to see that $F\in\Ld$ if and only if every entry $F_{ij}\in L^2(a,b)$. This induces a matrix-valued inner product for any two matrix-valued functions $F,G\in\Ld$, denoted by
\begin{equation}\label{inner}
   \langle F,G\rangle= \int_{a}^{b}F(x)G^*(x)dx.
\end{equation}
This is not an inner product in the common sense, but it has properties similar to the usual scalar inner products. It is also possible to define a scalar product between two matrix-valued functions (see \cite{DPS}), given by
$$
(F,G)=\mbox{Tr}\left(\langle F,G\rangle\right).
$$
Therefore, $\Ld$ with the norm $\|F\|=\mbox{Tr}\left(\langle F,F\rangle\right)^{1/2}$ is a Hilbert space and (\ref{inner}) is the inner product (in fact, the set of equivalence classes $F\sim G$ if $\|F-G\|=0$). Then we can module Fourier expansions of orthonormal systems in $\Ld$ (see pp. 7--8 in \cite{DPS}).

In a similar way, we can define the weighted spaces $\LWd$ of all matrix-valued functions in one variable with the inner product
\begin{equation}\label{innerW}
    \langle F,G\rangle_W=\int_{a}^{b}F(x)W(x)G^*(x)dx,
\end{equation}
where $dW$ is a weight matrix with a smooth density $W$ with respect to the Lebesgue measure, satisfying (1) $W(B)$ is positive semidefinite for any Borel set $B\in\mathbb{R}$, (2) $W$ has finite moments of every order and (3)
$\langle P,P\rangle_W$ is nonsingular if the leading coefficient of a matrix polynomial $P$ is nonsingular. Condition (3) above is necessary and sufficient to guarantee the existence of a sequence of matrix polynomials orthogonal with respect to \eqref{innerW} of degree $n$ with nonsingular leading coefficient, and it is fulfilled, in particular,
when $W$ is positive definite at an interval of the real line. We will assume, for simplicity, that there exist a real number $c\in(a,b)$ such that $W(c)=I$.

\bigskip

Now we introduce two matrices that are going to play a very important role in the rest of the paper. The first one is the $N\times N$ nilpotent matrix  (i.e. $A^N$=0)
\begin{equation}\label{AAA}
A=\sum_{j=1}^N\nu_jE_{j,j+1},\quad \nu_j\in\mathbb{R},
\end{equation}
where $E_{j,k}$ is a matrix with 1 at entry $(j,k)$ and 0
elsewhere, while the second is the diagonal matrix
\begin{equation}\label{JJ}
J=\sum_{j=1}^N(N-j)E_{j,j}.
\end{equation}
$A$ and $J$ satisfy the algebraic relation
\begin{equation*}\label{AlgRel}
    \mbox{ad}_AJ=[A,J]=-A,
\end{equation*}
where $\mbox{ad}^k, k\geq0$, denote the usual adjoint operators defined by $\mbox{ad}_X^0Y=Y$,
\begin{equation}\label{adj}
 \mbox{ad}_XY=XY-YX\quad\mbox{and}\quad\mbox{ad}_X^{k+1}Y=\mbox{ad}_X(\mbox{ad}_X^kY), \quad k\geq1.
\end{equation}
As a consequence
\begin{equation}\label{algrelk}
    \mbox{ad}_{A^k}J=-kA^k,\quad k=1,\ldots,N-1,
\end{equation}
as it is easy to see by induction.

\bigskip

We will now introduce real or complex valued functions with matrix argument. The Taylor series of an infinitely differentiable $f$ in a neighborhood of $x=0$ is
$$
f(x)=\sum_{j=0}^{\infty}f^{(j)}(0)\frac{x^j}{j!}.
$$
Whenever we write $f(X)$, for any $N\times N$ matrix $X$, we will mean the following matrix
$$
f(X)=\sum_{j=0}^{\infty}f^{(j)}(0)\frac{X^j}{j!}.
$$
It is clear that $f(X)g(X)=g(X)f(X)$ for any two real or complex valued functions $f,g$.

In particular we will be interested in the evaluation of $f(A)$ or $f(J)$, where $A$ and $J$ are defined in \eqref{AAA} and \eqref{JJ} respectively, for certain real or complex valued functions $f$. Observe first that, since $A$ is nilpotent, then
$$
f(A)=\sum_{j=0}^{N-1}f^{(j)}(0)\frac{A^j}{j!},
$$
i.e. $f(A)$ is a finite sum of linear combinations of the powers of $A$. It is clear that all real or complex valued algebraic or differential manipulations for functions $f(x)$ can be applied to the matrix $f(A)$. For instance, when we write $(I+A)^{-1}$, it will denote the matrix
$$
(I+A)^{-1}=\sum_{j=0}^{N-1}(-1)^jA^j.
$$

On the contrary, the diagonal matrix $f(J)$ will not be a finite sum of linear combinations of the powers of $J$. That is the case of the following diagonal matrices. For $f(x)=e^{i\frac{\pi}{2}kx}, k\in\mathbb{Z}$, the matrix $f(J)$ is defined by
\begin{equation}\label{iJJ}
    e^{i\frac{\pi}{2}kJ}=(i)^{kJ}=e^{kJ\log i}=\sum_{j=1}^N(i)^{k(N-j)}E_{j,j},\quad k\in\mathbb{Z},
\end{equation}
where $i$ is the imaginary unit. These matrices will play a very important role in the integral equations that we will study in Sections \ref{FIRST} and \ref{SECOND}. Observe that for $k=0$ in \eqref{iJJ} we have the identity matrix $I$. For $k=1$ we have the complex diagonal matrix
\begin{equation}\label{iJ}
    e^{i\frac{\pi}{2}J}=\sum_{j=1}^N(i)^{N-j}E_{j,j}=\begin{pmatrix}
  (i)^{N-1} &  &  & \\
   & \ddots &  &  \\
   &  & i &  \\
  & & & 1 \\
\end{pmatrix},
\end{equation}
and for $k=2$ we have
\begin{equation}\label{iJ2}
    e^{i\pi J}=\sum_{j=1}^N(-1)^{N-j}E_{j,j}=\begin{pmatrix}
  (-1)^{N-1} &  &  & \\
   & \ddots &  &  \\
   &  & -1 &  \\
  & & & 1 \\
\end{pmatrix},
\end{equation}
which is real and satisfies $e^{i\pi J}e^{i\pi J}=I$. For $k=3$ we have the inverse of \eqref{iJ}, while for $k\geq4$ all matrices will reduce to one of the matrices mentioned above. These diagonal matrices are going to play the same role that the imaginary unit $i$ plays in the scalar situation.

Observe that in \eqref{iJJ} we are taking the principal value of $\log$. However it is clear that if we take any other value we will have that $e^{i(\frac{\pi}{2}+2m\pi)kJ}=e^{i\frac{\pi}{2}kJ}$ for any $m\in\mathbb{Z}$, so it will not give any additional information.

Similarly we can define the following diagonal (and singular) matrices
\begin{equation*}\label{Sin}
    \sin\bigg(\frac{\pi}{2}kJ\bigg)=\frac{1}{2i}(e^{i\frac{\pi}{2}kJ}-e^{-i\frac{\pi}{2}kJ})=\sum_{j=1}^N\sin\bigg(\frac{\pi}{2}k(N-j)\bigg)E_{j,j},\quad k\in\mathbb{Z},
\end{equation*}
\begin{equation*}\label{Cos}
 \cos\bigg(\frac{\pi}{2}kJ\bigg)=\frac{1}{2}(e^{i\frac{\pi}{2}kJ}+e^{-i\frac{\pi}{2}kJ})=\sum_{j=1}^N\cos\bigg(\frac{\pi}{2}k(N-j)\bigg)E_{j,j} ,\quad k\in\mathbb{Z}.
\end{equation*}
In the special case of $k=1$ the explicit expression of $\sin\big(\frac{\pi}{2}J\big)$ and $\cos\big(\frac{\pi}{2}J\big)$ are the diagonal real matrices
\begin{equation}\label{Sink2}
    \sin\bigg(\frac{\pi}{2}J\bigg)=\begin{pmatrix}
                      \sin\big(\frac{\pi}{2}(N-1)\big) & &  & &  & & \\
                       &  \ddots &  & &  & & \\
                         &  & -1 &  \\
                           &  &  & 0 &  \\
                          &  &  &  & 1 &  \\
                          &  &  &  &  & 0 \\
                       \end{pmatrix},
\end{equation}
and
\begin{equation}\label{Cosk2}
\cos\bigg(\frac{\pi}{2}J\bigg)=\begin{pmatrix}
\cos\big(\frac{\pi}{2}(N-1)\big) & &  & &  & & \\
                        & \ddots &  & &  & & \\
                          &  & 0 &  \\
                          &  &  & -1 &  \\
                         &  &  &  & 0 &  \\
                         &  &  &  &  & 1 \\
                       \end{pmatrix}.
\end{equation}

This is the only significant case, since for $k=2$, $\sin(\pi J)=0$ and $\cos(\pi J)=e^{i\pi J}$ and for $k=3$, $\sin\big(\frac{3\pi}{2}J\big)=-\sin\big(\frac{\pi}{2}J\big)$ and $\cos\big(\frac{3\pi}{2}J\big)=\cos\big(\frac{\pi}{2}J\big)$.

We remark the following relations
\begin{eqnarray}\label{sincosrel}
\nonumber e^{i\frac{\pi}{2}J}\cos\bigg(\frac{\pi}{2}J\bigg)=\frac{1}{2}(I+e^{i\pi J}), & \quad& e^{i\frac{\pi}{2}J}\sin\bigg(\frac{\pi}{2}J\bigg)=\frac{1}{2i}(e^{i\pi J}-I), \\
  e^{i\pi J}\cos\bigg(\frac{\pi}{2}J\bigg)=\cos\bigg(\frac{\pi}{2}J\bigg), & \quad& e^{i\pi J}\sin\bigg(\frac{\pi}{2}J\bigg)=-\sin\bigg(\frac{\pi}{2}J\bigg), \\
 \nonumber e^{i\frac{3\pi}{2}J}\cos\bigg(\frac{\pi}{2}J\bigg)=\frac{1}{2}(I+e^{i\pi J}), &\quad & e^{i\frac{3\pi}{2}J}\sin\bigg(\frac{\pi}{2}J\bigg)=-\frac{1}{2i}(e^{i\pi J}-I).
\end{eqnarray}
We will use these matrices to derive real integral equations of matrix-valued orthogonal polynomials in Section \ref{FIRST}.

\bigskip

Now we will prove the following lemma, which gives commutativity relations between the powers of $A$ and $e^{i\frac{\pi}{2}kJ}$:
\begin{lema}\label{lema1}
The matrices $A$ and $e^{i\frac{\pi}{2}kJ}$, defined in \eqref{AAA} and \eqref{iJJ} respectively, satisfy the following algebraic relation
\begin{equation}\label{AiJ}
    e^{i\frac{\pi}{2}kJ}A^m=(i)^{km}A^me^{i\frac{\pi}{2}kJ},\quad m=1, 2,\ldots,N-1,\quad k\in\mathbb{Z}.
\end{equation}
\end{lema}
\pro
It is enough to prove this for $m=1$. For the rest of values of $m$ we apply recursively the formula for $m=1$. As a consequence of the definitions of $A$ and $e^{i\frac{\pi}{2}kJ}$ in \eqref{AAA} and \eqref{iJJ} respectively, we have
\begin{align*}
(i)^kAe^{i\frac{\pi}{2}kJ}=&(i)^k\bigg(\sum_{j=1}^{N-1}\nu_jE_{j,j+1}\bigg)\bigg(\sum_{l=1}^N(i)^{k(N-l)}E_{l,l}\bigg)= (i)^k\sum_{j=1}^{N-1}(i)^{k(N-j-1)}\nu_jE_{j,j+1}\\
=&\sum_{j=1}^{N-1}(i)^{k(N-j)}\nu_jE_{j,j+1}=\bigg(\sum_{j=1}^N(i)^{k(N-j)}E_{j,j}\bigg)\bigg(\sum_{l=1}^{N-1}\nu_lE_{l,l+1}\bigg)=e^{i\frac{\pi}{2}kJ}A,
\end{align*}
using that $E_{j,k}E_{h,m}=E_{j,m}$ for $k=h$ or 0 otherwise.
\hfill$\Box$

\section{Matrix-valued orthogonal functions}\label{MVOF}

As we mentioned in the Introduction, in the last few years many families of $N\times N$ matrix-valued orthogonal polynomials $(P_n)_n$ have been found along with their orthogonality measure $W$ satisfying second-order differential equations of the form
\begin{equation}\label{sode}
    P_n''(x)F_2(x)+P_n'(x)F_1(x)+P_n(x)F_0(x)=\Gamma_nP_n(x),
\end{equation}
where $F_2, F_1$ and $F_0$ are matrix polynomials (which do not depend on $n$) of degrees less than or equal to 2, 1 and 0, respectively, and $\Gamma_n$ are Hermitian matrices (if the family $(P_n)_n$ is orthonormal). These families are natural orthogonal systems in the weighted spaces $\LWd$.

Typically the weight matrices with this property can be factorized in the form
$$
W(x)=\rho(x)T(x)T^*(x),
$$
where $\rho$ is a scalar function (Hermite, Laguerre or Jacobi weight) and $T$ is a matrix-valued function which satisfies a first order differential equation with initial conditions of the form
\begin{equation}\label{GGG}
 T'(x)=G(x)T(x),\quad T(c)=I,
\end{equation}
for some $c\in(a,b)$. The differential coefficient $G$ in \eqref{GGG} is connected with the coefficients $F_2, F_1$ of the differential operator \eqref{sode} and $\rho$ (see the proof below of Theorem \ref{Phiteo} for more details or \cite{DG1}). Many examples have been found in the last years by solving a set of three symmetry differential equations that are equivalent to the second-order differential operator \eqref{sode} being self-adjoint with respect to the inner product \eqref{innerW} (see \cite{DG1, DG7, GPT5, D8, D3}). For a different approach, using matrix-valued spherical functions, see \cite{GPT6, PT1}.

One possibility for finding orthogonal systems in the space $\Ld$ is considering, for each family of matrix-valued orthogonal polynomials $(P_n)_n$ the following family
\begin{equation}\label{Phi}
 \Phi_n(x)=\rho^{1/2}(x)P_n(x)T(x),\quad n\geq0.
\end{equation}

Then $(\Phi_n)_n$ will be orthogonal with respect to the identity matrix $I$. Considering \eqref{sode} and \eqref{GGG} one verifies that the family $(\Phi_n)_n$ also satisfies very special differential equations. We will restrict ourselves to the special case where the differential coefficient $F_2(x)$ in \eqref{sode} is scalar, i.e. of the form $F_2(x)=f_2(x)I,$ for some real polynomial $f_2(x)$ of degree less than or equal to 2. This is the case considered in \cite{DG1}. Therefore we have the following

\begin{teor}\label{Phiteo}
The family $(\Phi_n)_n$ defined by \eqref{Phi}, satisfies the following second-order differential equation
\begin{equation}\label{Phisode}
    f_2(x)\Phi_n''(x)+f_2'(x)\Phi_n'(x)-\Phi_n(x)\bigg(\frac{1}{4}\bigg(\frac{(f_2(x)\rho'(x))'}{\rho(x)}+\bigg(\frac{f_2(x)\rho'(x)}{\rho(x)}\bigg)'\bigg)I+\chi(x)\bigg)=\Gamma_n\Phi_n(x),
\end{equation}
where $\chi(x)$ is the Hermitian matrix-valued function
\begin{equation*}\label{chi}
    \chi(x)=T^{-1}(x)\bigg(f_2(x)G'(x)+f_2(x)G^2(x)+\frac{(f_2(x)\rho(x))'}{\rho(x)}G(x)-F_0\bigg)T(x),
\end{equation*}
$G$ is the coefficient in \eqref{GGG} and $F_0$ and $\Gamma_n$ are the independent coefficient and the eigenvalue, respectively, of the differential equation \eqref{sode}.
\end{teor}
\pro
Differentiating twice \eqref{Phi} one gets expressions of $\rho^{1/2}P_n'T$ and $\rho^{1/2}P_n''T$ in terms of the derivatives of $\Phi_n$. Multiplying the differential equation \eqref{sode} on the right by $\rho^{1/2}T$ and substituting the formulas mentioned above we get \eqref{Phisode} after using $F_1=2f_2G+\frac{(f_2\rho)'}{\rho}$ (see formula (4.5) in \cite{DG1}). The matrix-valued function $\chi$ is Hermitian as a consequence of the third symmetry equation (see formula (4.12) in \cite{DG1} for details).
\hfill$\Box$
\bigskip

Observe that the coefficients of the differential operator \eqref{Phisode} are extremely simplified and are all scalar functions except for the independent coefficient or potential, which is Hermitian. In fact, in all the examples in the literature until now, the matrix $\chi$ is a diagonal matrix (see comment in pp. 93 of \cite{D3}).

The expression of the differential equation \eqref{Phisode} for each one the classical weights of Hermite, Laguerre and Jacobi is
\begin{enumerate}
  \item For $\rho(x)=e^{-x^2}$ and $f_2(x)=1$
\begin{equation}\label{Schrodinger}
    \Phi_n''(x)-\Phi_n(x)((x^2-1)I+\chi(x))=\Gamma_n\Phi_n(x).
\end{equation}
\item For $\rho(x)=x^{\alpha}e^{-x}, \alpha>-1,$ and $f_2(x)=x$
\begin{equation}\label{confluent}
x\Phi_n''(x)+\Phi_n'(x)-\Phi_n(x)\bigg(\bigg(\frac{x}{4}-\frac{\alpha+1}{2}+\frac{\alpha^2}{4x}\bigg)I+\chi(x)\bigg)=\Gamma_n\Phi_n(x).
\end{equation}
  \item For $\rho(x)=(1-x)^{\alpha}(1+x)^{\beta}, \alpha,\beta>-1,$ and $f_2(x)=1-x^2$
\begin{equation}\label{spheroidal}
(1-x^2)\Phi_n''(x)-2x\Phi_n'(x)-\Phi_n(x){\normalsize \mbox{$\bigg(\bigg(\frac{\alpha^2}{2(1-x)}+\frac{\beta^2}{2(1+x)}-\frac{(\alpha+\beta)(\alpha+\beta+2)}{2}\bigg)I+\chi(x)\bigg)$}}=\Gamma_n\Phi_n(x).
\end{equation}
\end{enumerate}
The explicit expression of the Hermitian matrix-valued function $\chi$ in each case can be found in Lemma 2.2 of \cite{DG7} for \eqref{Schrodinger} and \eqref{confluent} and in formula (2.5) of \cite{D3} for \eqref{spheroidal}.

Observe that the equation \eqref{Schrodinger} may be regarded of as a matrix-valued version of the one-dimensional \emph{Schr\"{o}dinger equation}. The equation \eqref{confluent} is related with a matrix-valued version of the \emph{confluent hypergeometric equation}, while the equation \eqref{spheroidal} is related with a matrix-valued version of the differential equation for \emph{spheroidal wave functions} (for special values of the parameters $\alpha$ and $\beta$). These differential equations in the scalar case are very important with diverse applications to physics, engineering and mathematical analysis itself (see \cite{CH, Le, B}).

\bigskip

In this paper we will focus only on families satisfying matrix-valued Schr\"{o}dinger differential equations as in \eqref{Schrodinger}, i.e. supported in the real line. For each of the families that we consider, we will prove that there exists an integral operator of the form \eqref{Inteqgen} with kernel $K(x,t)$ having the matrix-valued functions $(\Phi_n)_n$ as eigenfunctions. We will concentrate on Examples 5.1 and 5.2 of \cite{DG1}.

We will always denote by $(\Phi_n)_n$ the matrix-valued functions for all the examples, while there is no confusion about which family we are using from the context. In the case we are referring to different examples at the same time we will denote them by $(\Phi_{n,k})_n$, $k=1,2,$ where $(\Phi_{n,1})_n$ is the example in Section \ref{FIRST} and $(\Phi_{n,2})_n$ is the example in Section \ref{SECOND}.

\section{The first example}\label{FIRST}

Let $W$ be the following weight matrix
\begin{equation}\label{W1}
    W(x)=e^{-x^2}e^{Ax}e^{A^*x},\quad x\in\mathbb{R},
\end{equation}
where $A$ is the $N\times N$ nilpotent matrix \eqref{AAA}.

This example was considered for the first time in \cite{DG1}. The family of monic orthogonal polynomials $(\widehat{P}_n)_n$  satisfies a second-order differential equation as in \eqref{sode} with
$$
F_2(x)=I,\quad F_1(x)=-2xI+2A,\quad F_0(x)=A^2-2J,\quad \Gamma_n=-2nI+A^2-2J,
$$
where $J$ is the $N\times N$ diagonal matrix defined in \eqref{JJ}.

Now we construct a family of polynomials of the form $P_n(x)=L_n\widehat{P}_n(x)$ where the leading coefficient $L_n$ is chosen such that the eigenvalue $\Gamma_n$ transforms into one diagonal. A possible choice for this example is
\begin{equation*}\label{Ln1}
    L_n=e^{-A^2/4}.
\end{equation*}
It is straightforward to see, using the formula $e^{X}He^{-X}=\sum_{j\geq0}\mbox{ad}_X^j(H)/j!$ (see \eqref{adj} for definitions) and \eqref{algrelk} for $k=2$, that $L_n\Gamma_nL_n^{-1}=-2nI-2J$.

Hence, denoting
\begin{equation}\label{P1}
P_n(x)=e^{-A^2/4}\widehat{P}_n(x),
\end{equation}
the family
\begin{equation}\label{Phi1}
\Phi_n(x)=e^{-x^2/2}P_n(x)e^{Ax}
\end{equation}
satisfies a more convenient differential equation
\begin{equation}\label{Schrodinger1}
\Phi_n''(x)-\Phi_n(x)(x^2I+2J)+((2n+1)I+2J)\Phi_n(x)=0,
\end{equation}
as a consequence of Theorem \ref{Phiteo} (see \eqref{Schrodinger}). Observe that in this case we have $\chi(x)=2J$ (see Theorem 5.1 of \cite{DG1}). We remark that this differential equation is independent of the matrix $A$. Note that a similar differential equation (using the monic family) was derived in Section 6.2 of \cite{DG3}.

In order to prove the main result in this section (Theorem \ref{teor1} below) we need the following
\begin{lema}\label{lema2}
The following formula holds
\begin{equation}\label{formlema1}
    \frac{1}{\sqrt{2\pi}}\int_{-\infty}^{\infty}e^{-t^2/2}e^{-A^2/4}e^{At}e^{ixt}e^{i\frac{\pi}{2}J}dt=e^{i\frac{\pi}{2}J}e^{-x^2/2}e^{-A^2/4}e^{Ax},
\end{equation}
where $e^{i\frac{\pi}{2}J}$ is the diagonal matrix \eqref{iJ}.
\end{lema}
\pro
Expanding $e^{At}$ on the left hand side of the formula \eqref{formlema1}, using formula \eqref{AiJ} for $k=1$ and $e^{-A^2/4}e^{i\frac{\pi}{2}J}=e^{i\frac{\pi}{2}J}e^{A^2/4}$ (as a consequence of \eqref{AiJ}), and denoting $\widehat{H}_n(x)=(-1)^ne^{x^2/2}(e^{-x^2/2})^{(n)}$ the monic Hermite polynomials, we obtain
\begin{align*}
\frac{1}{\sqrt{2\pi}}\int_{-\infty}^{\infty}e^{-t^2/2}e^{-A^2/4}e^{At}e^{ixt}e^{i\frac{\pi}{2}J}dt=& \frac{1}{\sqrt{2\pi}}\sum_{j=0}^{N-1}e^{-A^2/4}A^je^{i\frac{\pi}{2}J}\int_{-\infty}^{\infty}e^{-t^2/2}\frac{t^j}{j!}e^{ixt}dt\\
=&\frac{1}{\sqrt{2\pi}}\sum_{j=0}^{N-1}e^{-A^2/4}e^{i\frac{\pi}{2}J}A^j(-i)^j\int_{-\infty}^{\infty}e^{-t^2/2}\frac{t^j}{j!}e^{ixt}dt\\
=&\sum_{j=0}^{N-1}e^{i\frac{\pi}{2}J}e^{A^2/4}\frac{A^j}{j!}(-i)^je^{-x^2/2}(i)^j\widehat{H}_j(x)\\
=&e^{i\frac{\pi}{2}J}e^{-x^2/2}e^{A^2/4}\sum_{j=0}^{N-1}\frac{A^j}{j!}\widehat{H}_j(x)\\
=&e^{i\frac{\pi}{2}J}e^{A^2/4}e^{-x^2/2}e^{Ax-A^2/2}=e^{i\frac{\pi}{2}J}e^{-x^2/2}e^{-A^2/4}e^{Ax},
\end{align*}
since the monic Hermite functions $\widehat{H}_n(x)e^{-x^2/2}$ are eigenfunctions of the Fourier transform with eigenvalue $(i)^n$ and the generating function for the monic Hermite polynomials $(\widehat{H}_n)_n$ is given by $\sum_{j=0}^{\infty}\widehat{H}_j(x)\frac{t^j}{j!}=e^{xt-t^2/2}$.
\hfill$\Box$
\bigskip

\begin{teor}\label{teor1}
The family of matrix-valued orthogonal functions $(\Phi_n)_n$ defined in \eqref{Phi1} satisfies the following integral equation
\begin{equation}\label{IntEq1}
\frac{1}{\sqrt{2\pi}}\int_{-\infty}^{\infty}\Phi_n(t)e^{ixt}e^{i\frac{\pi}{2}J}dt=(i)^ne^{i\frac{\pi}{2}J}\Phi_n(x),\quad x\in\mathbb{R},
\end{equation}
where $e^{i\frac{\pi}{2}J}$ is the diagonal matrix \eqref{iJ}.
\end{teor}

\pro Denote
\begin{equation}\label{FT1}
    \Psi_n(x)=\frac{1}{\sqrt{2\pi}}\int_{-\infty}^{\infty}\Phi_n(t)e^{ixt}e^{i\frac{\pi}{2}J}dt.
\end{equation}
By integration by parts using $\frac{d}{dx}e^{ixt}=ite^{ixt}$ and $\frac{d}{dt}e^{ixt}=ixe^{ixt}$ we get that $\Psi_n(x)$ satisfies the same second-order differential equation as $\Phi_n(x)$, i.e. \eqref{Schrodinger1}. Therefore $\Psi_n(x)$ can be written as $\Psi_n(x)=C_n\Phi_n(x)$, $n\geq0$, for some sequence of nonsingular diagonal matrices $C_n$. This follows if we expand $\Psi_n(x)=\sum_{k=0}^{\infty}C_{n,k}\Phi_k(x)$. Since $\Psi_n(x)$ satisfies \eqref{Schrodinger1} and $(\Phi_n)_n$ is a system of linearly independent matrix-valued functions this is equivalent to say that $(n-k)C_{n,k}+JC_{n,k}-C_{n,k}J=0$ for all $n,k\geq0$. But now it is straightforward to conclude that if $n\neq k$ then $C_{n,k}=0$ and if $n=k$ then $C_{n,n}\doteq C_n$ must be a diagonal matrix, since $J$ has simple spectrum.

Now we will use Lemma \ref{lema2}, which is exactly the case $n=0$ in \eqref{FT1}. After differentiating the right hand side of the formula \eqref{formlema1} $n$ times with respect to $x$ one obtains
\begin{equation}\label{expr1}
    e^{i\frac{\pi}{2}J}e^{-A^2/4}\frac{d^n}{dx^n}(e^{-x^2/2}e^{Ax})=\frac{(i)^ne^{-A^2/4}}{\sqrt{2\pi}}\int_{-\infty}^{\infty}e^{-t^2/2}t^ne^{At}e^{ixt}e^{i\frac{\pi}{2}J}dt.
\end{equation}
Using the Leibniz's formula, the left hand side of the formula \eqref{expr1} can be also written as
\begin{equation}\label{expr2}
    e^{i\frac{\pi}{2}J}e^{-A^2/4}\frac{d^n}{dx^n}(e^{-x^2/2}e^{Ax})=e^{i\frac{\pi}{2}J}e^{-A^2/4}((-1)^nx^nI+\cdots)e^{-x^2/2}e^{Ax}.
\end{equation}
Writing $\Phi_n(x)=e^{-x^2/2}P_n(x)e^{Ax}=e^{-x^2/2}e^{-A^2/4}(x^nI+\cdots)e^{Ax}$, by linearity, and using the right hand side of \eqref{expr1} and \eqref{expr2} one obtains
\begin{align*}
\Psi_n(x)=\frac{1}{\sqrt{2\pi}}\int_{-\infty}^{\infty}\Phi_n(t)e^{ixt}e^{i\frac{\pi}{2}J}dt&=\frac{e^{-A^2/4}}{\sqrt{2\pi}}\int_{-\infty}^{\infty}(t^nI+\cdots)e^{At}e^{ixt}e^{i\frac{\pi}{2}J}dt\\
&=(i)^ne^{i\frac{\pi}{2}J}e^{-A^2/4}(x^nI+\cdots)e^{-x^2/2}e^{Ax}.
\end{align*}
Since we already know that $\Psi_n(x)=C_n\Phi_n(x)$, equating the leading degree of $x$ in the formula above we obtain $C_n=(i)^ne^{i\frac{\pi}{2}J}$.
\hfill$\Box$

\bigskip

Observe that the integral equation \eqref{IntEq1} is also independent of the matrix $A$.

The integral equation will have important consequences:

\begin{coro}\label{coro11}
The family of matrix-valued orthogonal functions $(\Phi_n)_n$ satisfies the following symmetry condition
\begin{equation}\label{symmeqphi}
    \Phi_n(x)=(-1)^ne^{i\pi J}\Phi_n(-x)e^{i\pi J},
\end{equation}
where $e^{i\pi J}$ is the real diagonal matrix \eqref{iJ2}. Consequently, the family of matrix-valued orthogonal polynomials $P_n(x)=e^{-A^2/4}\widehat{P}_n(x), n\geq0,$ satisfies the same symmetry condition, i.e.
\begin{equation}\label{symmeqpn}
    P_n(x)=(-1)^ne^{i\pi J}P_n(-x)e^{i\pi J}.
\end{equation}
\end{coro}
\pro
From \eqref{IntEq1} multiplying on the left by the eigenvalue $(i)^ne^{i\frac{\pi}{2}J}$ and substituting again by the same formula we get
\begin{align*}
(-1)^ne^{i\pi J}\Phi_n(x)&=\frac{1}{\sqrt{2\pi}}\int_{-\infty}^{\infty}(i)^ne^{i\frac{\pi}{2}J}\Phi_n(t)e^{ixt}e^{i\frac{\pi}{2}J}dt =\frac{1}{2\pi}\int_{-\infty}^{\infty}\int_{-\infty}^{\infty}\Phi_n(z)e^{izt}e^{ixt}e^{i\pi J}dzdt\\
&=\int_{-\infty}^{\infty}\Phi_n(z)\delta(x+z)e^{i\pi J}dz=\Phi_n(-x)e^{i\pi J},
\end{align*}
using standard Fourier analysis, where $\delta$ is the standard Dirac delta function. Therefore \eqref{symmeqphi} holds. The formula \eqref{symmeqpn} holds from the observation that $e^{-Ax}e^{i\pi J}e^{-Ax}=e^{i\pi J}$, as a consequence of \eqref{AiJ}.
\hfill$\Box$

\bigskip

\begin{coro}\label{coro12}
The family of matrix-valued orthogonal polynomials $P_n(x)=e^{-A^2/4}\widehat{P}_n(x)$, $n\geq0,$ with respect to the weight matrix \eqref{W1} satisfies the following integral equation
\begin{equation}\label{IntEqPn}
e^{-x^2/2}P_n(x)e^{Ax}=\frac{(-i)^n}{\sqrt{2\pi}}e^{-i\frac{\pi}{2}J}\int_{-\infty}^{\infty}P_n(t)e^{-t^2/2}e^{ixt}e^{At}e^{i\frac{\pi}{2}J}dt,
\end{equation}
where $e^{i\frac{\pi}{2}J}$ is the diagonal matrix \eqref{iJ}.
\end{coro}
\pro
Immediate from \eqref{IntEq1} using \eqref{Phi1}.
\hfill$\Box$

\bigskip

We remark here that the family $(P_n)_n$ is a solution of the integral equation \eqref{Inteqgen} with kernel $K(x,t)=e^{(x+it)^2}e^{At}e^{i\frac{\pi}{2}J}e^{-Ax}$.

\bigskip

Observe that \eqref{IntEqPn} is an integral equation which depends on complex values. In the scalar case it is well known that the Hermite polynomials $(H_n)_n$ satisfy real integral equations in terms of the kernels $\cos(xt)$ and $\sin(xt)$ (see for instance \cite{Le}). This is possible since the Hermite polynomials satisfy the symmetry condition $H_n(x)=(-1)^nH_n(-x)$. In our case, we have a different symmetry condition \eqref{symmeqpn}, so it will not follow the same lines as in the scalar case. However, it is possible to derive real integral equations for $(P_n)_n$:

\begin{coro}\label{coro15}
The family of matrix-valued orthogonal polynomials $P_n(x)=e^{-A^2/4}\widehat{P}_n(x)$, $n\geq0,$ with respect to the weight matrix \eqref{W1} satisfies the following real integral equations
\begin{equation}\label{IntEqPneven}
   e^{-x^2/2}(e^{i\pi J}\pm I)P_{n}(x)e^{Ax}C_{\pm}=\frac{(-1)^{\lfloor\frac{n}{2}\rfloor}}{\sqrt{2\pi}}C_{\pm}\int_{-\infty}^{\infty}e^{-t^2/2}k_n(x,t)P_{n}(t) e^{At}dt(e^{i\pi J}\pm I),
\end{equation}
and
\begin{equation}\label{IntEqPnodd}
     e^{-x^2/2}(e^{i\pi J}\pm I)P_{n}(x)e^{Ax}C_{\mp}=\frac{\pm (-1)^{\lfloor\frac{n}{2}\rfloor}}{\sqrt{2\pi}}C_{\pm}\int_{-\infty}^{\infty}e^{-t^2/2}k_{n+1}(x,t)P_{n}(t)e^{At}dt (e^{i\pi J}\mp I),
\end{equation}
where $e^{i\pi J}$ is defined in \eqref{iJ2}, $C_{-}=\sin\big(\frac{\pi}{2}J\big)$, $C_{+}=\cos\big(\frac{\pi}{2}J\big)$  and are both defined in \eqref{Sink2} and \eqref{Cosk2} respectively, $\lfloor \cdot\rfloor$ denotes the floor function and
\begin{equation*}
    k_n(x,t)=\left\{\begin{array}{cc}
               \cos(xt), & \mbox{if $n$ is even}, \\
               \sin(xt), & \mbox{if $n$ is odd}.
             \end{array}\right.
\end{equation*}
\end{coro}
\pro
From the symmetry condition \eqref{symmeqpn} and using formula \eqref{AiJ} we have that $$e^{-x^2/2}e^{i\frac{\pi}{2}J}P_n(x)e^{Ax}=(-1)^ne^{i\pi J}[e^{-x^2/2}e^{-i\frac{\pi}{2}J}P_n(-x)e^{-Ax}]e^{i\pi J}.$$ Therefore, the evaluation of \eqref{IntEqPn} at $x$ and $-x$ and using the previous formula gives
\begin{align*}
2e^{-x^2/2}e^{i\frac{\pi}{2}J}P_n(x)e^{Ax}&=\frac{(i)^n}{\sqrt{2\pi}}\bigg[\int_{-\infty}^{\infty}e^{-t^2/2}\cos{xt}\bigg((-1)^nP_n(t)e^{At}e^{i\frac{\pi}{2}J}+e^{i\pi J}P_n(t)e^{At}e^{i\frac{\pi}{2}J}e^{i\pi J}\bigg)dt\\
&+i\int_{-\infty}^{\infty}e^{-t^2/2}\sin{xt}\bigg((-1)^nP_n(t)e^{At}e^{i\frac{\pi}{2}J}-e^{i\pi J}P_n(t)e^{At}e^{i\frac{\pi}{2}J}e^{i\pi J}\bigg)dt\bigg].
\end{align*}
Now we multiply on the right and on the left by appropriate matrices such that the elements in the big parenthesis of the formula above are equal and the other vanishes for even or odd values of $n$. This will become clear below. These matrices are $\sin\big(\frac{\pi}{2}J\big)$ and $\cos\big(\frac{\pi}{2}J\big)$, defined in \eqref{Sink2} and \eqref{Cosk2} respectively. There are four possible combinations of multiplying these matrices on the right and on the left. This is why we will get eventually \emph{eight} formulas, considering even and odd values of $n$. We will show the case when we multiply on the left by $\cos\big(\frac{\pi}{2}J\big)$ and on the right by $\sin\big(\frac{\pi}{2}J\big)$. From the analysis below the rest of formulas can be derived in a similar way.

From the previous formula and using relations \eqref{sincosrel} one gets
\begin{align*}
e^{-x^2/2}(e^{i\pi J}+I)P_n(x)e^{Ax}\sin\bigg(\frac{\pi}{2}J\bigg)&=\frac{(i)^n}{2i\sqrt{2\pi}}\cos\bigg(\frac{\pi}{2}J\bigg)\bigg[\int_{-\infty}^{\infty}e^{-t^2/2}\cos(xt)((-1)^n-1)P_n(t)e^{At}dt\\
&+i\int_{-\infty}^{\infty}e^{-t^2/2}\sin(xt)((-1)^n+1)P_n(t)e^{At}dt\bigg](e^{i\pi J}-I).
\end{align*}
Now it is clear that for even or odd values of $n$, either the integral with $\cos(xt)$ or $\sin(xt)$ will vanish and that these integral equations are real. Therefore
\begin{equation*}
    e^{-x^2/2}(e^{i\pi J}+I)P_{2n}(x)e^{Ax}\sin\bigg(\frac{\pi}{2}J\bigg)=\frac{(-1)^n}{\sqrt{2\pi}}\cos\bigg(\frac{\pi}{2}J\bigg)\int_{-\infty}^{\infty}e^{-t^2/2}\sin(xt)P_{2n}(t)e^{At}dt(e^{i\pi J}-I),
\end{equation*}
\begin{equation*}
    e^{-x^2/2}(e^{i\pi J}+I)P_{2n+1}(x)e^{Ax}\sin\bigg(\frac{\pi}{2}J\bigg)=\frac{(-1)^{n+1}}{\sqrt{2\pi}}\cos\bigg(\frac{\pi}{2}J\bigg)\int_{-\infty}^{\infty}e^{-t^2/2}\cos(xt)P_{2n+1}(t)e^{At}dt(e^{i\pi J}-I),
\end{equation*}
which are formulas \eqref{IntEqPnodd} for the positive sign and both odd and even values of $n$.
\hfill$\Box$
\bigskip

Observe that \eqref{IntEqPneven} and \eqref{IntEqPnodd} give \emph{eight} real integral equations of the family $(P_n)_n$ according to positive or negative sign and odd or even values of $n$. The reason for eight formulas is the structure of the matrices $\sin\big(\frac{\pi}{2}J\big)$ and $\cos\big(\frac{\pi}{2}J\big)$, defined in \eqref{Sink2} and \eqref{Cosk2} respectively. These are diagonal singular matrices and whenever we multiply them on the right or on the left we only get a description of some of the entries of $P_n$ (which is a full matrix in general). It is exactly this combination of multiplying on the right and on the left by these matrices that allows us to have a formula for every entry of $P_n$. In other words, if we multiply on the left by $\cos\big(\frac{\pi}{2}J\big)$ and on the right by both $\cos\big(\frac{\pi}{2}J\big)$ or $\sin\big(\frac{\pi}{2}J\big)$, we get a description of the $N, N-2, N-4,\ldots$ rows of $P_n$ while if we multiply on the left by $\sin\big(\frac{\pi}{2}J\big)$ and on the right by both $\cos\big(\frac{\pi}{2}J\big)$ or $\sin\big(\frac{\pi}{2}J\big)$, we get a description of the $N-1, N-3, N-5,\ldots$ rows of $P_n$, respectively. Therefore all rows of $P_n$ are covered.

\begin{nota}
We have found that the family of polynomials $(P_n)_n$ defined in \eqref{P1} (and the family we will study in Section \ref{SECOND} defined in \eqref{P2}) simplifies considerably many structural formulas. For instance, the norms of $(P_n)_n$ (and consequently the norms of $(\Phi_n)_n$ with respect to the inner product \eqref{inner}) are diagonal. Also the coefficients of the three-term recurrence relation are considerably simplified. A detailed study of these and other structural formulas for $(P_n)_n$ will be discussed in future publications.
\end{nota}

\subsection{A detailed study of the case $N=2$}\label{GRAPHS1}

In this section we will study in detail the properly normalized matrix-valued functions $(\Phi_n)_n$ for the special case of $N=2$. We will show their relationship with the scalar Hermite or wave functions, derive some structural formulas and plot graphs of the diagonal entries of $(\Phi_n\Phi_n^*)_n$ for some values of $n$ and $\nu_1$ (the only free parameter in the matrix $A$, see \eqref{AAA}).

\bigskip

Using Theorem 5.1 of \cite{DG3} one can derive an explicit expression for the family of matrix-valued orthogonal polynomials $(P_n)_n$ defined in \eqref{P1}. They can be written as
\begin{equation*}\label{PnN2}
    P_n(x)=\frac{1}{2^n}\begin{pmatrix}
                          H_n(x) & -n\nu_1H_{n-1}(x) \\
                          -\frac{n\nu_1}{\gamma_n}H_{n-1}(x)  & \frac{1}{\gamma_n}(H_n(x)+n\nu_1^2xH_{n-1}(x)) \\
                        \end{pmatrix},
\end{equation*}
where $H_n(x)=(-1)^ne^{x^2}(e^{-x^2})^{(n)}$, $n\geq0,$ are the classical Hermite polynomials, and $(\gamma_n)$ is a sequence of real positive numbers defined by
\begin{equation}\label{gamma1}
    \gamma_n=1+\frac{n}{2}\nu_1^2,\quad n\geq0.
\end{equation}

We can normalize this family since the norms of $(P_n)_n$ are diagonal
$$
\|P_n\|_W^2=\frac{n!\sqrt{\pi}}{2^n}\begin{pmatrix}
                                    \gamma_{n+1} & 0 \\
                                    0 & 1/\gamma_n \\
                                  \end{pmatrix}.
$$
Therefore the normalized matrix-valued functions $\widetilde{\Phi}_n=\|P_n\|_W^{-1}\Phi_n$ can be written as
\begin{equation}\label{PhinN2}
    \widetilde{\Phi}_n(x)=\begin{pmatrix}
                          \psi_n(x)/\sqrt{\gamma_{n+1}} & \nu_1\sqrt{\frac{n+1}{2\gamma_{n+1}}}\psi_{n+1}(x) \\
                          -\nu_1\sqrt{\frac{n}{2\gamma_n}}\psi_{n-1}(x)  & \psi_n(x)/\sqrt{\gamma_{n}}\\
                        \end{pmatrix},
\end{equation}
where $\psi_n(x)=\frac{1}{\sqrt{2^nn!\sqrt{\pi}}}e^{-x^2/2}H_n(x)$ are the normalized Hermite or wave functions. Hence
\begin{equation}\label{Phin2N2}
    \widetilde{\Phi}_n(x)\widetilde{\Phi}_n^*(x)=\begin{pmatrix}
                          \psi_{n+1}^2(x)+\frac{1}{\gamma_{n+1}}(\psi_{n}^2(x)-\psi_{n+1}^2(x)) & \frac{\nu_1}{\sqrt{\gamma_{n}\gamma_{n+1}}}\psi_{n}(x)\psi_{n}'(x) \\
                          \frac{\nu_1}{\sqrt{\gamma_{n}\gamma_{n+1}}}\psi_{n}(x)\psi_{n}'(x)  & \psi_{n-1}^2(x)+\frac{1}{\gamma_{n}}(\psi_{n}^2(x)-\psi_{n-1}^2(x))\\
                        \end{pmatrix},\quad n\geq0.
\end{equation}
In the expressions above we are using standard properties of Hermite polynomials and wave functions like $H_{n+1}(x)=2xH_n(x)-2nH_{n-1}(x)$ and $\psi_{n}'(x)=\sqrt{\frac{n+1}{2}}\psi_{n+1}(x)-\sqrt{\frac{n}{2}}\psi_{n-1}(x)$. Observe from the definition of $(\gamma_n)_n$ in \eqref{gamma1} that the diagonal entries of \eqref{Phin2N2} are probability densities depending on one free parameter $\nu_1$.

The effect of the parameter $\nu_1$ can be studied from the explicit expression of \eqref{PhinN2} and \eqref{Phin2N2}. For instance, as $\nu_1$ tends to 0 we observe that the diagonal entries of \eqref{PhinN2} converge to two copies of the Hermite functions while the off diagonal entries tend to 0. If $\nu_1$ is large, the diagonal entries get small compared with the off diagonal entries (for the first values of $n$).

In Figure 1 and 2 we have plotted the diagonal entries of $\widetilde{\Phi}_n(x)\widetilde{\Phi}_n^*(x)$, which are probability densities in $\RR$, for every $n=0,1,\ldots,5$, and for the special choice of $\nu_1=1$. The numbers on the graphics correspond to the value of $n$. It is well known that in the scalar situation there exist $n$ points (the zeros of the Hermite polynomials) where the probability density $\psi_n(x)\psi_n^*(x)$ is exactly 0. In our case our two probability densities never vanish.

\begin{figure}[h]
\begin{minipage}{7.5cm}
\includegraphics[width=6.5cm,height=7.63cm,angle=-90,scale=1.05]{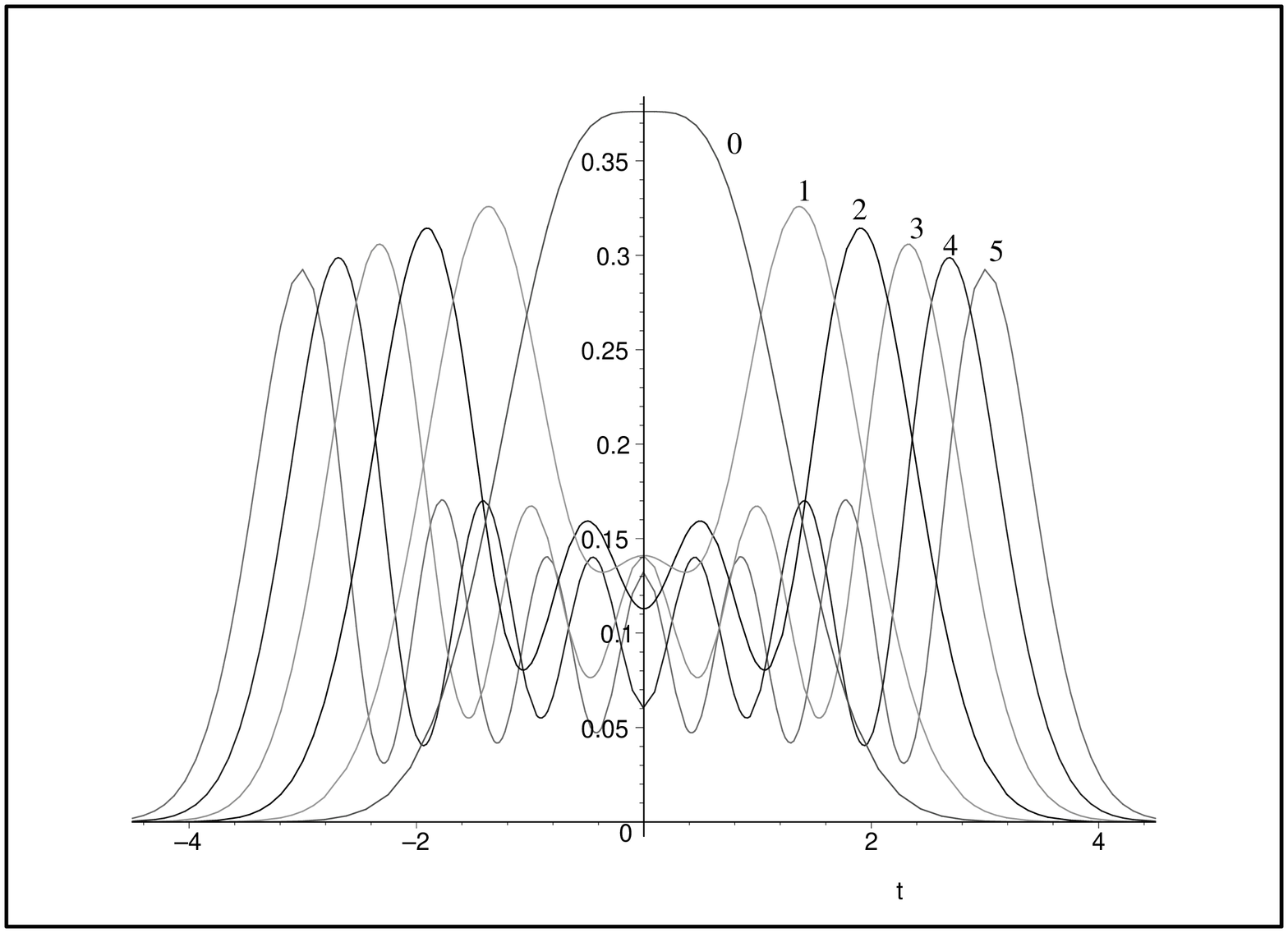}
\caption{$(\widetilde{\Phi}_{n}\widetilde{\Phi}_{n}^*)_{11}, n=0,\ldots,5, \nu_1=1$}
\end{minipage}
\    \ \
\begin{minipage}{7.5cm}
\includegraphics[width=6.5cm,height=7.63cm,angle=-90,scale=1.05]{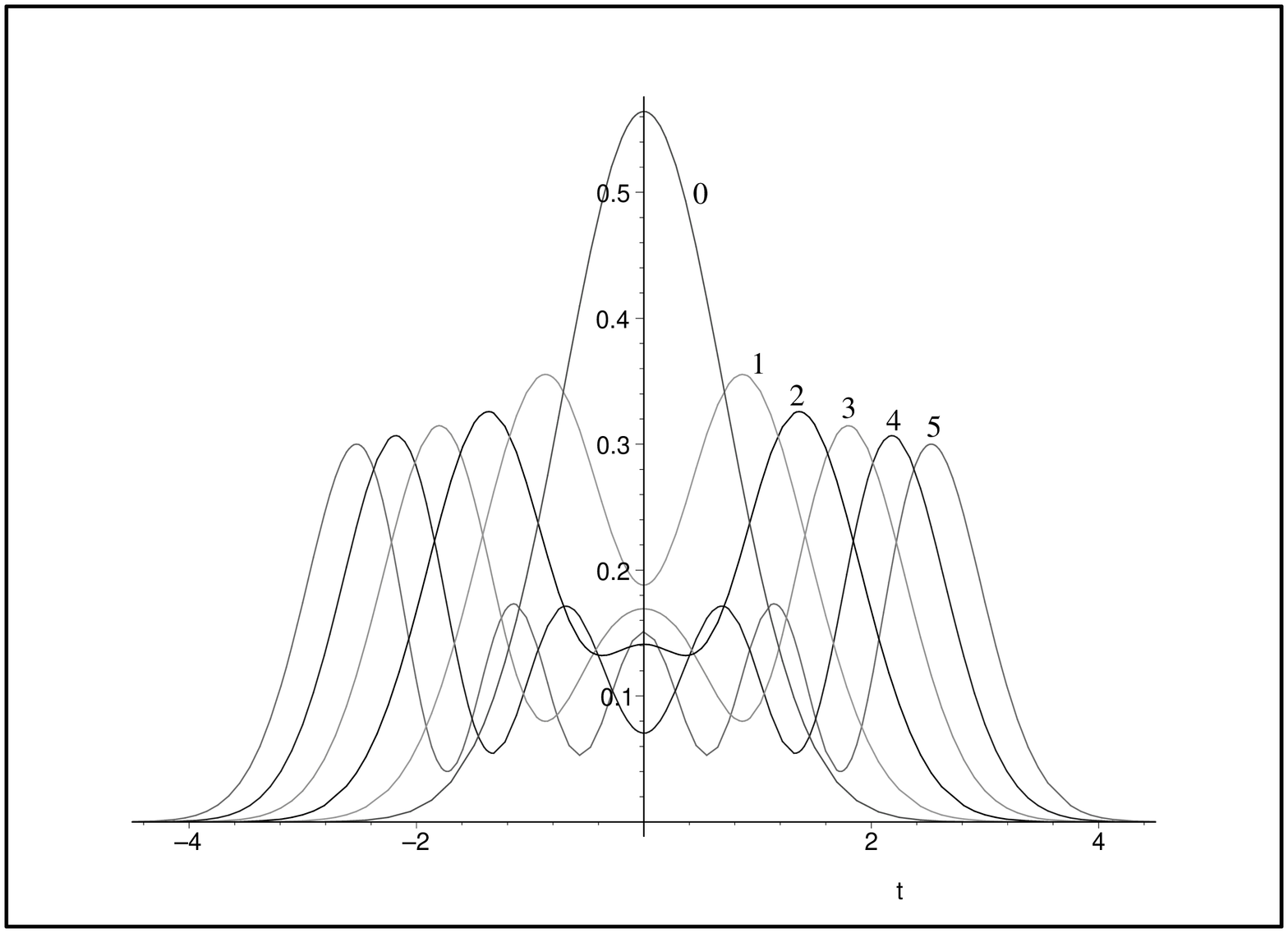}
\caption{$(\widetilde{\Phi}_{n}\widetilde{\Phi}_{n}^*)_{22}, n=0,\ldots,5, \nu_1=1$}
\end{minipage}
\end{figure}


\bigskip

Finally we compute the following formulas
\begin{equation}\label{Molec}
    (x^kI)_{nm}\doteq\int_{-\infty}^{\infty}x^k\widetilde{\Phi}_n(x)\widetilde{\Phi}_m^*(x)dx,\quad n,m\geq0,\quad k=1,2.
\end{equation}
In the scalar situation, we have, as a consequence of the three-term recurrence relation among Hermite polynomials, that
\begin{equation}\label{waveform1}
(x)_{nm}=\sqrt{\frac{n}{2}}\delta_{m,n-1}+\sqrt{\frac{n+1}{2}}\delta_{m,n+1},
\end{equation}
\begin{equation}\label{waveform2}
(x^2)_{nm}=\frac{1}{2}\sqrt{n(n-1)}\delta_{m,n-2}+(n+1/2)\delta_{m,n}+\frac{1}{2}\sqrt{(n+1)(n+2)}\delta_{m,n+2}.
\end{equation}

We can use the formulas in the scalar case to calculate the matrix-valued formulas \eqref{Molec}. From \eqref{PhinN2} we can obtain an explicit expression of $\widetilde{\Phi}_n\widetilde{\Phi}_m^*$:
\begin{equation*}\label{Phinm2N2}
    \widetilde{\Phi}_n\widetilde{\Phi}_m^*=\begin{pmatrix}
                          \frac{\psi_n\psi_m}{\sqrt{\gamma_{n+1}\gamma_{m+1}}}+\frac{\nu_1^2}{2}\sqrt{\frac{(n+1)(m+1)}{\gamma_{n+1}\gamma_{m+1}}} \psi_{n+1}\psi_{m+1} & \frac{\nu_1}{\sqrt{\gamma_{m}\gamma_{n+1}}}(\sqrt{\frac{n+1}{2}}\psi_{n+1}\psi_{m}-\sqrt{\frac{m}{2}}\psi_{n}\psi_{m-1}) \\
 \frac{\nu_1}{\sqrt{\gamma_{n}\gamma_{m+1}}}(\sqrt{\frac{m+1}{2}}\psi_{m+1}\psi_{n}-\sqrt{\frac{n}{2}}\psi_{m}\psi_{n-1})  & \frac{\psi_n\psi_m}{\sqrt{\gamma_{n}\gamma_{m}}}+\frac{\nu_1^2}{2}\sqrt{\frac{nm}{\gamma_{n}\gamma_{m}}} \psi_{n-1}\psi_{m-1}\\
                        \end{pmatrix}.
\end{equation*}
And using formulas \eqref{waveform1} and \eqref{waveform2} for the scalar wave functions we obtain
\begin{align*}\label{waveMform1}
(xI)_{nm}=&\begin{pmatrix}
                                                                         \sqrt{\frac{n\gamma_{n+1}}{2\gamma_n}} & 0 \\
                                                                         0 & \sqrt{\frac{n\gamma_{n-1}}{2\gamma_n}} \\
                                                                       \end{pmatrix}
\delta_{m,n-1}+\begin{pmatrix}
                 0 & \frac{\nu_1}{2\sqrt{\gamma_n\gamma_{n+1}}} \\
                 \frac{\nu_1}{2\sqrt{\gamma_n\gamma_{n+1}}}  & 0 \\
               \end{pmatrix}
\delta_{m,n}
\\
&
+\begin{pmatrix}
                                                                         \sqrt{\frac{(n+1)\gamma_{n+2}}{2\gamma_{n+1}}} & 0 \\
                                                                         0 & \sqrt{\frac{(n+1)\gamma_{n}}{2\gamma_{n+1}}} \\
                                                                       \end{pmatrix}\delta_{m,n+1},
\end{align*}
and
\begin{align*}
(x^2I)_{nm}&=\begin{pmatrix}
                                                                         \frac{1}{2}\sqrt{\frac{n(n-1)\gamma_{n+1}}{\gamma_{n-1}}} & 0 \\
                                                                         0 & \frac{1}{2}\sqrt{\frac{n(n-1)\gamma_{n-2}}{\gamma_{n}}} \\
                                                                       \end{pmatrix}\delta_{m,n-2}+\begin{pmatrix}
                 0 & \nu_1\sqrt{\frac{n}{2\gamma_{n-1}\gamma_{n+1}}} \\
                 \frac{\nu_1}{\gamma_{n}}\sqrt{\frac{n}{2}}  & 0 \\
               \end{pmatrix}\delta_{m,n-1}\\
&+\begin{pmatrix}
    n+\frac{3}{2}-\frac{1}{\gamma_{n+1}} & 0 \\
    0 & n-\frac{1}{2}+\frac{1}{\gamma_{n}} \\
  \end{pmatrix}
\delta_{m,n}+\begin{pmatrix}
                 0 & \frac{\nu_1}{\gamma_{n+1}}\sqrt{\frac{n+1}{2}} \\
                 \nu_1\sqrt{\frac{n+1}{2\gamma_{n}\gamma_{n+2}}}  & 0 \\
               \end{pmatrix}\delta_{m,n+1}\\
&+\begin{pmatrix}
                                                                         \frac{1}{2}\sqrt{\frac{(n+1)(n+2)\gamma_{n+3}}{\gamma_{n+1}}} & 0 \\
                                                                         0 & \frac{1}{2}\sqrt{\frac{(n+1)(n+2)\gamma_{n}}{\gamma_{n+2}}} \\
                                                                       \end{pmatrix}\delta_{m,n+2},
\end{align*}
where $(\gamma_n)_n$ is the sequence \eqref{gamma1}. Observe that for $k=1$ we get an extra term when $n=m$, something that did not happen in the scalar situation. Likewise for $k=2$ we get two new extra terms.

In the scalar situation the tridiagonal matrix $(x)$ given in \eqref{waveform1} can be  viewed as the matrix of the homomorphism $f\mapsto xf$ in $L^2(\RR)$ with respect to the basis $(\psi_n)_n$. In the matrix case for this example we have that the block tridiagonal matrix $(xI)$ is the matrix of the homomorphism $F\mapsto xF$ in $\LdR$ with respect to the basis $(\widetilde{\Phi}_n)_n$. In this case, $(xI)$ is a 4 diagonal semi-infinite matrix with zeros in the main diagonal of the form
\begin{equation*}\label{xOp1}
    (xI)=\begin{pmatrix}
      0 & \star & \star &  &  &  &  \\
      \star & 0 &0 & \star&  &  &  \\
      \star & 0 & 0& \star & \star&  &  \\
       & \star& \star & 0 & 0 &\star &  \\
 & & \star& 0 & 0 & \star &\star &  \\
      & &  & \ddots & \ddots & \ddots & \ddots & \ddots \\
    \end{pmatrix},
\end{equation*}
where a $\star$ means a nonzero entry.

\section{The second example}\label{SECOND}

Let $W$ be the following weight matrix
\begin{equation}\label{W2}
    W(x)=e^{-x^2}e^{A(I+A)^{-1}x^2}e^{(I+A^*)^{-1}A^*x^2},\quad x\in\mathbb{R},
\end{equation}
where $A$ is the $N\times N$ nilpotent matrix \eqref{AAA}. This weight matrix was considered for the first time in Example 5.2 of \cite{DG1}, where the authors use the notation $B=A(I+A)^{-1}$. Although this new notation looks more complicated it will simplify considerably all computations in this section. All proofs follow the same lines as in Section \ref{FIRST}.

The family of monic orthogonal polynomials $(\widehat{P}_n)_n$  satisfies a second-order differential equation as in \eqref{sode} with
$$
F_2(x)=I,\quad F_1(x)=2(2A(I+A)^{-1}-I)x,\quad F_0(x)=2A(I+A)^{-1}-4J,$$
$$
\Gamma_n=2(2A(I+A)^{-1}-I)n+2A(I+A)^{-1}-4J,
$$
where $J$ is the $N\times N$ diagonal matrix defined in \eqref{JJ}.

Now we construct a family of polynomials of the form $P_n(x)=L_n\widehat{P}_n(x)$ where the leading coefficient $L_n$ is chosen such that the eigenvalue $\Gamma_n$ transforms into one diagonal. Observe that the eigenvalue has $n$ dependence off the main diagonal. Therefore the coefficient $L_n$ will depend on $n$. A natural candidate in this case is
\begin{equation*}\label{Ln2}
    L_n=[(I+A)^{-1/2}]^{2n+1},
\end{equation*}
as a consequence the algebraic relation $\log(I+A)J-J\log(I+A)=-A(I+A)^{-1}$ which can be proved using \eqref{algrelk} (just expanding in power series).

Hence, denoting
\begin{equation}\label{P2}
P_n(x)=[(I+A)^{-1/2}]^{2n+1}\widehat{P}_n(x),
\end{equation}
the family
\begin{equation}\label{Phi12}
\Phi_n(x)=e^{-x^2/2}P_n(x)e^{A(I+A)^{-1}x^2}
\end{equation}
satisfies a more convenient differential equation
\begin{equation*}\label{Schrodinger12}
\Phi_n''(x)-\Phi_n(x)(x^2I+4J)+((2n+1)I+4J)\Phi_n(x)=0.
\end{equation*}
as a consequence of Theorem \ref{Phiteo} (see \eqref{Schrodinger}). Observe that in this case we have $\chi(x)=4J$ (see Theorem 5.2 of \cite{DG1}). We remark again that this differential equation is independent of the matrix $A$. Note that a similar differential equation (using the monic family) was derived in Section 6.2 of \cite{DG3}.

\bigskip

As before, in order to prove the main result in this section, we need the following
\begin{lema}\label{lema22}
The following formula holds
\begin{equation*}\label{formlema12}
    \frac{1}{\sqrt{2\pi}}\int_{-\infty}^{\infty}e^{-t^2/2}(I+A)^{-1/2}e^{A(I+A)^{-1}t^2}e^{ixt}e^{i\pi J}dt=e^{i\pi J}e^{-x^2/2}(I+A)^{-1/2}e^{A(I+A)^{-1}x^2},
\end{equation*}
where $e^{i\pi J}$ is defined in \eqref{iJ2}.
\end{lema}
\pro
Although the proof is more elaborated than the one in Lemma \ref{lema2}, it follows the same procedure using standard real analysis of power series of functions and the formulas $e^{i\pi J}A^k[(I+A)^{-1}]^k=(-1)^kA^k[(I+A)^{-1}]^ke^{i\pi J}$ and $(I+A)^{1/2}e^{i\pi J}=e^{i\pi J}(I-A)^{1/2}$, which hold using Lemma \ref{lema1}.
\hfill$\Box$

Therefore we have the following
\begin{teor}\label{teor2}
The family of matrix-valued orthogonal functions $(\Phi_n)_n$ defined in \eqref{Phi12} satisfies the following integral equation
\begin{equation*}\label{IntEq12}
\frac{1}{\sqrt{2\pi}}\int_{-\infty}^{\infty}\Phi_n(t)e^{ixt}e^{i\pi J}dt=(i)^ne^{i\pi J}\Phi_n(x),
\end{equation*}
where $e^{i\pi J}$ is defined in \eqref{iJ2}.
\end{teor}
\pro
The only difference with respect to the proof of Theorem \ref{teor1} is that the leading coefficient of $P_n$ depends on $n$ and that the formula \eqref{expr2} is now
\begin{align*}
\frac{d^n}{dx^n}(e^{-x^2/2}e^{A(I+A)^{-1}x^2})&=((-1)^n[2A(I+A)^{-1}-I]^nx^nI+\cdots)e^{-x^2/2}e^{A(I+A)^{-1}x^2})\\
&=((-1)^n[(I-A)(I+A)^{-1}]^nx^nI+\cdots)e^{-x^2/2}e^{A(I+A)^{-1}x^2}).
\end{align*}
The rest follows the same arguments as in Theorem \ref{teor1}.
\hfill$\Box$

\bigskip

In a similar way we have the following corollaries, from which we omit the proofs since they are exactly the same as in the example in Section \ref{FIRST}.

\begin{coro}\label{coro112}
The family of matrix-valued orthogonal functions $(\Phi_n)_n$ and the family of matrix-valued orthogonal polynomials $(P_n)_n$ satisfy the following symmetric conditions
\begin{equation}\label{symmeqphi2}
    \Phi_n(x)=(-1)^n\Phi_n(-x),\quad\quad P_n(x)=(-1)^nP_n(-x).
\end{equation}
\end{coro}
Observe now that the symmetry conditions are exactly the same as the classical Hermite polynomials.

\begin{coro}\label{coro122}
The family of matrix-valued orthogonal polynomials $(P_n)_n$, $n\geq0,$ with respect to the weight matrix \eqref{W2}, satisfies the following integral equation
\begin{equation*}\label{IntEqPn2}
e^{-x^2/2}P_n(x)e^{A(I+A)^{-1}x^2}=\frac{(-i)^n}{\sqrt{2\pi}}e^{i\pi J}\int_{-\infty}^{\infty}P_n(t)e^{-t^2/2}e^{ixt}e^{A(I+A)^{-1}t^2}e^{i\pi J}dt,
\end{equation*}
where $e^{i\pi J}$ is the diagonal matrix \eqref{iJ2}. Moreover, we have the following real integral equations
\begin{equation*}\label{IntEqPneven2}
   e^{-x^2/2}e^{i\pi J}P_{2n}(x)e^{A(I+A)^{-1}x^2}=\frac{(-1)^n}{\sqrt{2\pi}}\int_{-\infty}^{\infty}e^{-t^2/2}P_{2n}(t)e^{A(I+A)^{-1}t^2}\cos(xt)e^{i\pi J} dt,
\end{equation*}
\begin{equation*}\label{IntEqPnodd2}
     e^{-x^2/2}e^{i\pi J}P_{2n+1}(x)e^{A(I+A)^{-1}x^2}=\frac{(-1)^n}{\sqrt{2\pi}}\int_{-\infty}^{\infty}e^{-t^2/2}P_{2n+1}(t)e^{A(I+A)^{-1}t^2}\sin(xt) e^{i\pi J}dt.
\end{equation*}
\end{coro}

Observe now, because of \eqref{symmeqphi2}, that there are only two real integral equations, unlike the example studied in Section \ref{FIRST}.

\subsection{A detailed study of the case $N=2$}\label{GRAPHS2}

Using Theorem 5.1 of \cite{DG3} one can derive an explicit expression for the family of matrix-valued orthogonal polynomials $(P_n)_n$ defined in \eqref{P2}. They can be written as
\begin{equation*}\label{PnN22}
    P_n(x)=\frac{1}{2^n}\begin{pmatrix}
                          H_n(x) & -\nu_1((n+1/2)H_n(x)+n(n-1)H_{n-2}(x)) \\
                          -n(n-1)\nu_1H_{n-2}(x)/\gamma_n & H_n(x)/\gamma_n+n(n-1)\nu_1^2x^2H_{n-2}(x)/\gamma_n \\
                        \end{pmatrix},
\end{equation*}
where $H_n(x)=(-1)^ne^{x^2}(e^{-x^2})^{(n)}$, $n\geq0$ are the classical Hermite polynomials, and $(\gamma_n)_n$ is a sequence of real positive numbers defined by
\begin{equation}\label{gamma2}
    \gamma_n=1+\frac{\nu_1^2}{2}\binom{n}{2},\quad n\geq0.
\end{equation}
We can normalize this family since the norms of $(P_n)_n$ are diagonal
$$
\|P_n\|_W^2=\frac{n!\sqrt{\pi}}{2^n}\begin{pmatrix}
                                    \gamma_{n+2} & 0 \\
                                    0 & 1/\gamma_n \\
                                  \end{pmatrix}.
$$
Therefore the normalized matrix-valued functions $\widetilde{\Phi}_n=\|P_n\|_W^{-1}\Phi_n$ can be written in the following form
\begin{equation}\label{PhinN22}
    \widetilde{\Phi}_n(x)=\begin{pmatrix}
                          \psi_n(x)/\sqrt{\gamma_{n+2}} & \frac{\nu_1}{2}\sqrt{\frac{(n+1)(n+2)}{\gamma_{n+2}}}\psi_{n+2}(x) \\
                          -\frac{\nu_1}{2}\sqrt{\frac{n(n-1)}{\gamma_n}}\psi_{n-2}(x)  & \psi_n(x)/\sqrt{\gamma_{n}}\\
                        \end{pmatrix},
\end{equation}
where $\psi_n(x)=\frac{1}{\sqrt{2^nn!\sqrt{\pi}}}e^{-x^2/2}H_n(x)$ are the normalized Hermite functions. Therefore
\begin{equation}\label{Phin2N22}
    \widetilde{\Phi}_n(x)\widetilde{\Phi}_n^*(x)=\begin{pmatrix}
                        \psi_{n+2}^2(x)+\frac{1}{\gamma_{n+2}}(\psi_n^2(x)-\psi_{n+2}^2(x)) & -\frac{\nu_1}{2\sqrt{\gamma_{n}\gamma_{n+2}}}\psi_{n}(x)(\psi_{n}(x)+2x\psi_{n}'(x)) \\
                          -\frac{\nu_1}{2\sqrt{\gamma_{n}\gamma_{n+2}}}\psi_{n}(x)(\psi_{n}(x)+2x\psi_{n}'(x))  & \psi_{n-2}^2(x)+\frac{1}{\gamma_{n}}(\psi_n^2(x)-\psi_{n-2}^2(x))\\
                        \end{pmatrix}.
\end{equation}
Observe from the definition of $(\gamma_n)_n$ in \eqref{gamma2} that the diagonal entries of \eqref{Phin2N22} are probability densities depending on one free parameter $\nu_1$. The effect of the parameter $\nu_1$ is similar to the example studied in Section \ref{FIRST}.

In Figure 3 and 4 we have plotted the diagonal entries of $\widetilde{\Phi}_n(x)\widetilde{\Phi}_n^*(x)$ for $n=0,1,\ldots,5$, and for the special choice of $\nu_1=1/2$. The numbers on the graphics correspond to the value of $n$.

\begin{figure}[h]
\begin{minipage}{7.5cm}
\includegraphics[width=6.5cm,height=7.63cm,angle=-90,scale=1.05]{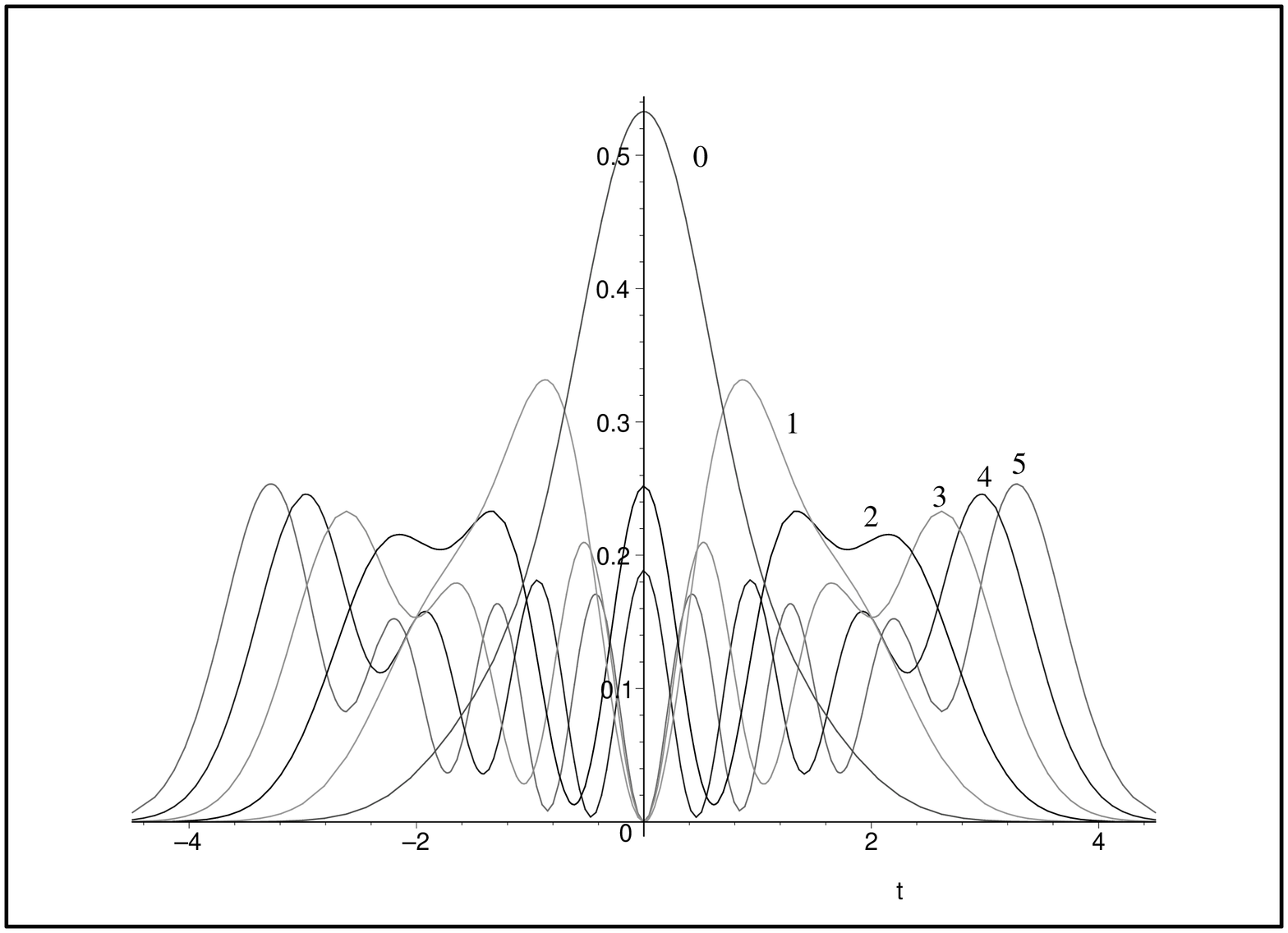}
\caption{$(\widetilde{\Phi}_{n}\widetilde{\Phi}_{n}^*)_{11}, n=0,...,5, \nu_1=1/2$}
\end{minipage}
\    \ \
\begin{minipage}{7.5cm}
\includegraphics[width=6.5cm,height=7.63cm,angle=-90,scale=1.05]{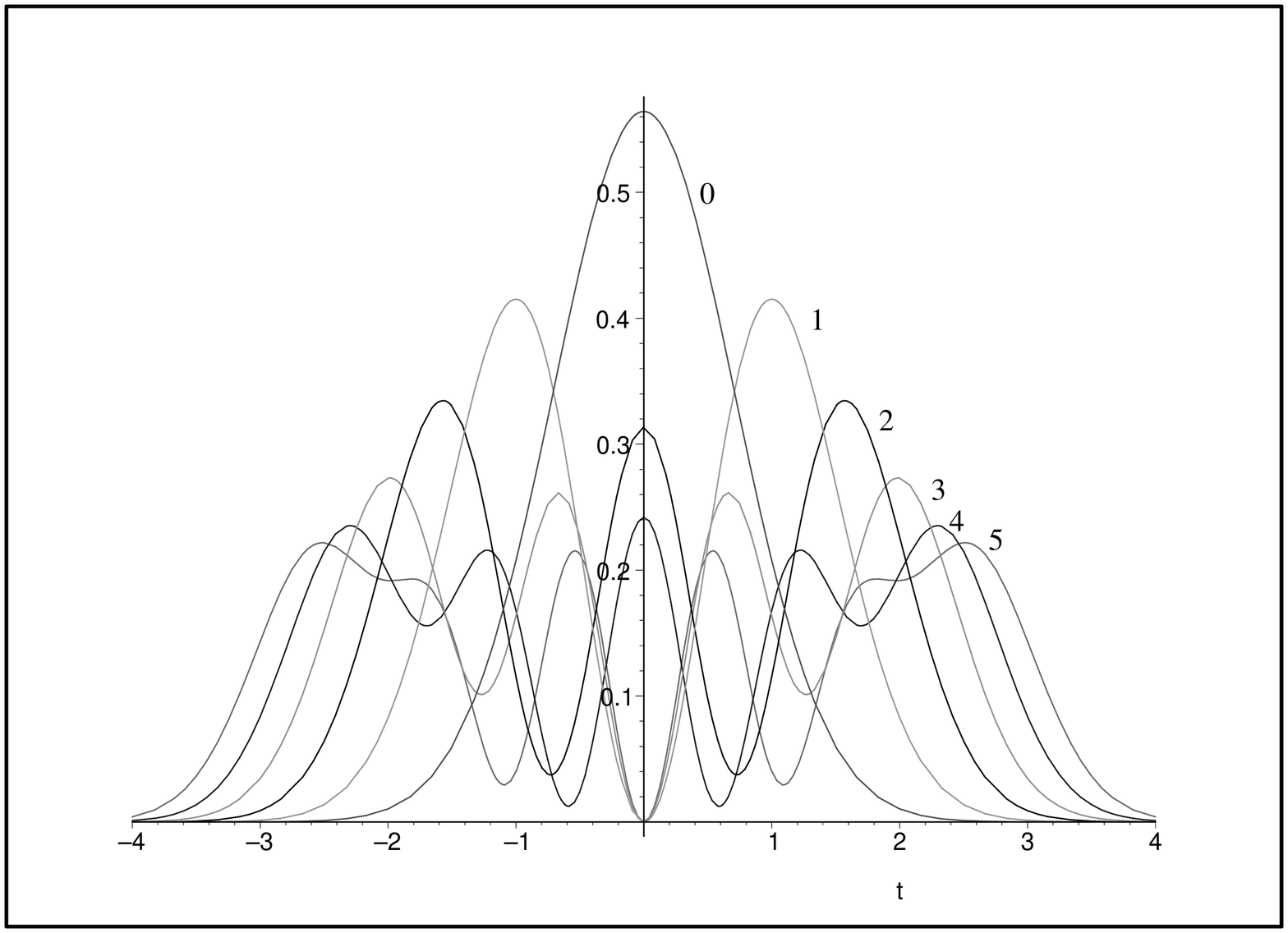}
\caption{$(\widetilde{\Phi}_{n}\widetilde{\Phi}_{n}^*)_{22}, n=0,...,5, \nu_1=1/2$}
\end{minipage}
\end{figure}

\bigskip

Finally we compute the following formulas
\begin{equation*}\label{Molec2}
    (x^kI)_{nm}\doteq\int_{-\infty}^{\infty}x^k\widetilde{\Phi}_n(x)\widetilde{\Phi}_m^*(x)dx,\quad n,m\geq0,\quad k=1,2,
\end{equation*}
in a similar way as in the previous example using formulas \eqref{waveform1} and
\eqref{waveform2}.

From \eqref{PhinN22} we can obtain an explicit expression of $\widetilde{\Phi}_n\widetilde{\Phi}_m^*$:
\begin{align*}
\widetilde{\Phi}_n(x)\widetilde{\Phi}_m^*(x)=&\left(\begin{array}{cc}
                                          \frac{\psi_n(x)\psi_m(x)}{\sqrt{\gamma_{n+2}\gamma_{m+2}}}+\frac{\nu_1^2}{4}\sqrt{\frac{(n+1)(n+2)(m+1)(m+2)}{\gamma_{n+2}\gamma_{m+2}}} \psi_{n+2}(x)\psi_{m+2}(x) & * \\
                                           \frac{\nu_1}{\sqrt{\gamma_{n}\gamma_{m+2}}}\bigg(\frac{\sqrt{(m+1)(m+2)}}{2}\psi_{m+2}(x)\psi_{n}(x)-\frac{\sqrt{n(n-1)}}{2}\psi_{m}(x)\psi_{n-2}(x)\bigg)& * \\
                                        \end{array}\right.\\
& \left.\begin{array}{cc}
                                          * &\frac{\nu_1}{\sqrt{\gamma_{m}\gamma_{n+2}}}\bigg(\frac{\sqrt{(n+1)(n+2)}}{2}\psi_{n+2}(x)\psi_{m}(x)-\frac{\sqrt{m(m-1)}}{2}\psi_{n}(x)\psi_{m-2}(x)\bigg)  \\
                                          * & \frac{\psi_n(x)\psi_m(x)}{\sqrt{\gamma_{n}\gamma_{m}}}+\frac{\nu_1^2}{4}\sqrt{\frac{n(n-1)m(m-1)}{\gamma_{n}\gamma_{m}}} \psi_{n-2}(x)\psi_{m-2}(x)\\
                                        \end{array}\right).
\end{align*}
Therefore we obtain
\begin{align*}
(xI)_{nm}=&\begin{pmatrix}
                                                                         \sqrt{\frac{n\gamma_{n+2}}{2\gamma_{n+1}}} & 0 \\
                                                                         \nu_1\sqrt{\frac{n}{2\gamma_{n}\gamma_{n+1}}} & \sqrt{\frac{n\gamma_{n-1}}{2\gamma_n}} \\
                                                                       \end{pmatrix}
\delta_{m,n-1}+\begin{pmatrix}
                                                                         \sqrt{\frac{(n+1)\gamma_{n+3}}{2\gamma_{n+2}}} & \nu_1\sqrt{\frac{n+1}{2\gamma_{n+1}\gamma_{n+2}}} \\
                                                                         0 & \sqrt{\frac{(n+1)\gamma_{n}}{2\gamma_{n+1}}} \\
                                                                       \end{pmatrix}\delta_{m,n+1},
\end{align*}
and
\begin{align*}
(x^2I)_{nm}&=\begin{pmatrix}
                                                                         \frac{1}{2}\sqrt{\frac{n(n-1)\gamma_{n+2}}{\gamma_{n}}} & 0 \\
                                                                         \frac{\nu_1\sqrt{n(n-1)}}{\gamma_n} & \frac{1}{2}\sqrt{\frac{n(n-1)\gamma_{n-2}}{\gamma_{n}}} \\
                                                                       \end{pmatrix}\delta_{m,n-2}+\begin{pmatrix}
    n+\frac{5}{2}-\frac{2}{\gamma_{n+2}} & \frac{\nu_1(2n+1)}{2\sqrt{\gamma_n\gamma_{n+2}}} \\
    \frac{\nu_1(2n+1)}{2\sqrt{\gamma_n\gamma_{n+2}}} & n-\frac{3}{2}+\frac{2}{\gamma_{n}} \\
  \end{pmatrix}
\delta_{m,n}\\&+\begin{pmatrix}
                                                                         \frac{1}{2}\sqrt{\frac{(n+1)(n+2)\gamma_{n+4}}{\gamma_{n+2}}} & \frac{\nu_1\sqrt{(n+1)(n+2)}}{\gamma_{n+2}} \\
                                                                         0 & \frac{1}{2}\sqrt{\frac{(n+1)(n+2)\gamma_{n}}{\gamma_{n+2}}} \\
                                                                       \end{pmatrix}\delta_{m,n+2},
\end{align*}
where $(\gamma_n)_n$ is the sequence \eqref{gamma2}.

\bigskip

Now the block tridiagonal matrix $(xI)$ is the matrix of the homomorphism $F\mapsto xF$ in $\LdR$ with respect to the basis $(\widetilde{\Phi}_n)_n$. In this case, $(xI)$ is a 7 diagonal semi-infinite matrix with zeros in the first three main diagonals of the form
\begin{equation*}\label{xOp2}
    (xI)=\begin{pmatrix}
      0 & 0 & \star & \star &  &  &  \\
      0 & 0 &0 & \star& 0 &  &  \\
      \star & 0 & 0& 0& \star& \star &  \\
      \star & \star& 0 & 0 & 0 &\star & 0 \\
 & 0& \star& 0 & 0 & 0 &\star & \star \\
      & & \ddots & \ddots & \ddots & \ddots & \ddots & \ddots & \ddots\\
    \end{pmatrix}.
\end{equation*}
where a $\star$ means a nonzero entry.

\section{Concluding remarks}\label{CON}

In this paper we have shown the first examples of integral operators having families of matrix-valued orthogonal functions as eigenfunctions, which at the same time are eigenfunctions of a second-order differential operator of Schr\"odinger type. We have focused on two examples supported in the real line appeared in \cite{DG1}. Certainly other examples supported in the real line have been found as well. For instance, the example in Theorem 1.1 in \cite{D3} and the example included in the Appendix A.1 in \cite{D3}. These examples are generalizations of the example studied in Section \ref{FIRST}. We have found that these examples also satisfy similar integral equations of the form \eqref{Inteqgen}, along with similar structural formulas and integral equations for the matrix-valued orthogonal polynomials.

\bigskip

These integral operators are slight modifications of the usual Fourier transform, but not necessarily the scalar Fourier transform. Instead, we can perform Fourier analysis for matrix-valued functions using the families of matrix-valued orthogonal functions introduced in Sections \ref{FIRST} and \ref{SECOND}. This is different from the usual Fourier analysis where the kernel is $K(x,t)=e^{ixt}I$, and the transform is applied entry by entry. For instance, for the examples $\Phi_{n,k}$, $k=1,2$, introduced in Sections \ref{FIRST} and \ref{SECOND} respectively we can recover any matrix-valued function $F\in\LdR$ that can be written in the form
\begin{equation*}\label{FF}
    F(x)=\sum_{n=0}^{\infty}C_{n,k}\Phi_{n,k}(x),\quad C_{n,k}=\langle F,\Phi_{n,k}\rangle,\quad k=1,2,
\end{equation*}
with the inner product $\langle \cdot,\cdot\rangle$ defined in \eqref{inner}. This representation is in fact always possible for any $F\in\LdR$, since our families of eigenfunctions $(\Phi_{n,k})_n$, $k=1,2$,  are \emph{complete} in $\LdR$. We omit the details, but this is a consequence of having our weight matrices \eqref{W1} and \eqref{W2} as nonsingular matrices that die out exponentially at infinity to the zero matrix, that is
$$
W(x)=O(e^{-\alpha |x|})\quad\mbox{for some}\quad\alpha>0\quad\mbox{as}\quad |x|\rightarrow\infty.
$$

Hence it is possible to extend to the matrix case all the results necessary to proof the completeness of Hermite (and Laguerre) polynomials (see, for instance, the Theorem appeared in \cite{He} or Section 6.5 of \cite{AAR}) and use Fourier expansions (see pp. 7--8 in \cite{DPS}) to conclude that our families $(\Phi_{n,k})_n$, $k=1,2$, are complete in $\LdR$.

Therefore, define the following two integral transforms of $F$
\begin{equation*}\label{FouTran}
    (F \mathcal{F}_k)(x)=\frac{1}{\sqrt{2\pi}}\int_{-\infty}^{\infty}F(t)e^{ixt}e^{i\frac{\pi}{2}kJ}dt,\quad k=1,2,
\end{equation*}
where $e^{i\frac{\pi}{2}kJ}$ is defined in \eqref{iJJ} and $\mathcal{F}_k$ acting on the right means that the kernel is multiplied on the right. The case $k=0$ is the usual Fourier transform entry by entry but if we multiply on the right by $e^{i\frac{\pi}{2}kJ}$ we can recover $F$ in a nicer way using the eigenfunctions $(\Phi_{n,k})_n$, $k=1,2$. If we define the inverse of the previous two integral transforms as
\begin{equation*}\label{InvFouTran}
    (F \mathcal{F}^{-1}_k)(x)=\frac{1}{\sqrt{2\pi}}\int_{-\infty}^{\infty}F(t)e^{-ixt}e^{-i\frac{\pi}{2}kJ}dt,\quad k=1,2,
\end{equation*}
it is easy to see that $(\Phi_{n,k}\mathcal{F}_k)(x)=(i)^ne^{i\frac{\pi}{2}kJ}\Phi_{n,k}(x)$ and $(\Phi_{n,k}\mathcal{F}^{-1}_k)(x)=(-i)^ne^{i\frac{\pi}{2}kJ}\Phi_{n,k}(x),$ using the symmetry condition $\Phi_{n,k}(x)=(-1)^ne^{i\pi kJ}\Phi_{n,k}(-x)e^{i\pi kJ}$ (see \eqref{symmeqphi} and \eqref{symmeqphi2}). Therefore
\begin{equation*}\label{recF}
    F(x)=((F\mathcal{F}_k)\mathcal{F}^{-1}_k)(x).
\end{equation*}

That means that Fourier analysis of matrix-valued functions in $\LdR$ can be studied, at least, in two different ways according to the Fourier type transform that we consider and its corresponding eigenfunctions, not only applying Fourier transform entry by entry.

%
%

\end{document}